\documentclass{amsart}

\usepackage{amssymb}

\usepackage[all]{xy}

\newtheorem{thm}{Theorem}%[section]
\newtheorem{thm-defn}{Theorem-Definition}%[section]

\newtheorem{cor}[thm]{Corollary}

\newtheorem{prop}[thm]{Proposition}

\newtheorem{principle}[thm]{Principle} %!!!!!!!!!!!!!!!!!!!!!!
\theoremstyle{definition}
\newtheorem{defn}[thm]{Definition}
\newtheorem{prel-defn}[thm]{Preliminary Definition}

\newtheorem{say}[thm]{}
\newtheorem{exmp}[thm]{Example}

\newtheorem{aside}[thm]{Aside}          
\newtheorem*{ack}{Acknowledgments}      % \renewcommand{\theack}{} 

\newtheorem{defn-thm}[thm]{Definition--Theorem}  %!!!!!!!!!!!!!!!!!!!!!!!!
\newtheorem{defn-lem}[thm]{Definition--Lemma}  %!!!!!!!!!!!!!!!!!!!!!!!!
\theoremstyle{remark}

\renewcommand{\c}[0]{{\mathbb C}}  

\renewcommand{\o}[0]{{\mathcal O}} 
\newcommand{\z}[0]{{\mathbb Z}}

\renewcommand{\a}[0]{{\mathbb A}}

\newcommand{\p}[0]{{\mathbb P}}
\newcommand{\f}[0]{{\mathbb F}}
\newcommand{\q}[0]{{\mathbb Q}}
\newcommand{\map}[0]{\dasharrow}
\newcommand{\qtq}[1]{\quad\mbox{#1}\quad}
\newcommand{\spec}[0]{\operatorname{Spec}}

\newcommand{\supp}[0]{\operatorname{Supp}}    
    
\newcommand{\codim}[0]{\operatorname{codim}}    
    
\newcommand{\proj}[0]{\operatorname{Proj}}

\newcommand{\coker}[0]{\operatorname{coker}}

\newcommand{\aut}[0]{\operatorname{Aut}}

\newcommand{\tors}[0]{\operatorname{tors}}

\newcommand{\hilb}[0]{\operatorname{Hilb}}
\newcommand{\univ}[0]{\operatorname{Univ}}

\newcommand{\rdown}[1]{\lfloor{#1}\rfloor}

\newcommand{\simq}[0]{\sim_{\q}}

\newcommand{\mor}[0]{\operatorname{Mor}} 
 
\newcommand{\tsum}[0]{\textstyle{\sum}}

\def\into{\DOTSB\lhook\joinrel\to}

\def\loccoh#1.#2.#3.#4.{H^{#1}_{#2}(#3,#4)}

\DeclareMathAlphabet{\mathchanc}{OT1}{pzc}%
                                {m}{it}

\begin{document}
\bibliographystyle{amsalpha}

\today

\title{Moduli of varieties of general type}
\author{J\'anos Koll\'ar}

\maketitle
\tableofcontents

\section{Introduction}

\begin{say}[Moduli theories]\label{basic.mod.say}
The development of a good moduli theory consists of four basic steps.
\medskip

(\ref{basic.mod.say}.1) {\it Identify a class of  objects 
whose moduli theory is nice.} 

In some cases the answer is obvious, for instance, we should study the
moduli theory of smooth projective  curves.
In other cases, it took some time to understand what the correct 
objects should be: Abelian varieties should be replaced by
{\em polarized}  Abelian varieties, K3 surfaces by
{\em marked} K3 surfaces and vector bundles by
{\em stable} vector bundles. We see later that
smooth, projective varieties of general type should be replaced
by their canonical models. We discuss this in Section \ref{sec.2}.

In all these cases, we get non-proper moduli spaces.
 This is inconvenient for many reasons,
for instance
it is hard to count objects with various properties. 
To remedy this, one should look for a compactification whose
points have clear geometric meaning. 

\medskip
(\ref{basic.mod.say}.2)  {\it Choose a larger class of objects 
 to  form a proper moduli space.}

The choice of these objects is usually neither obvious nor unique.
It was not until the 1960's that the importance of this step was
understood and 
 stable curves and semi-stable sheaves
were identified and studied in detail.
For surfaces of general type the right class was described in \cite{ksb}
and for polarized  Abelian varieties in \cite{ale-abvar}.
The solution for varieties of general type is treated in  Section \ref{sec.3}.
\medskip

(\ref{basic.mod.say}.3) {\it Establish the correct moduli functor.}

Once a class of objects ${\mathbf V}$ is established, the 
corresponding  moduli functor is ususally declared to consist of all
flat families whose fibers are in ${\mathbf V}$. 
However, for  varieties of general type, allowing all
flat families gives the {\em wrong}
 moduli functor. The problem
 and the solution are analyzed in Section \ref{sec.4}.

\medskip

(\ref{basic.mod.say}.4) {\it Study the resulting moduli spaces and
their applications.}

The moduli of curves and the moduli of semi-stable vector bundles on curves
appear in many contexts and by now established themselves as one of the
richest  applications of algebraic geometry.
We are only at the beginning of the development of the moduli of
higher dimensional varieties; the basic results are outlined in
Section \ref{sec.5}.  I hope to see many more
applications in the future.
\end{say}

These notes are intended to give a survey of the subject,
stressing key examples and results. The forthcoming book
\cite{k-m-book} aims to give a complete treatment.

\begin{defn}[Moduli functors] \label{fine.mod.def.intro}
Let ${\mathbf V}$ be a ``reasonable'' class of projective varieties
(or schemes, or sheaves, or ...).
For the next definition we only need to  assume that if 
 $K\supset k$ is a field extension then 
a $k$-variety $X_k$ is in ${\mathbf V}$ iff
 $X_K:=X_k\times_{\spec k}\spec K$ is  in ${\mathbf V}$.
Define the 
corresponding moduli functor  as
$$
Varieties_{\mathbf V}(T):=
\left\{
\begin{array}{c}
\mbox{Flat families $X\to T$ such that}\\
\mbox{every fiber is in ${\mathbf V}$,}\\
\mbox{modulo isomorphisms over $T$.}
\end{array}
\right\}\
\eqno{(\ref{fine.mod.def.intro}.1)}
$$
(As noted in (\ref{basic.mod.say}.3), we will need to impose additional 
restrictions eventually, but for now let us ignore these.)
\end{defn}

\begin{say}[Moduli spaces]
\label{course.mod.varss.defn.intro}

We say that a 
scheme ${\rm Moduli}_{\mathbf V}$, or, more precisely, a
flat morphism 
$$
u:\univ_{\mathbf V}\to {\rm Moduli}_{\mathbf V}
$$
is a {\it fine moduli space} for the functor $Varieties_{\mathbf V}$
if, for every scheme $T$, pulling back gives an equality
$$
Varieties_{\mathbf V}(T)=\mor\bigl(T, {\rm Moduli}_{\mathbf V}\bigr).
%\eqno{(\ref{fine.mod.def.intro}.2)}
$$
Our aim is to understand all families whose fibers are in ${\mathbf V}$
and a fine moduli space presents the answer in the most succinct way.

Applying the definition
to $T=\spec K$, where $K$ is a field,
we see that every fiber of 
$u:\univ_{\mathbf V}\to {\rm Moduli}_{\mathbf V}$
is in ${\mathbf V}$
and the $K$-points of 
${\rm Moduli}_{\mathbf V}$ are in
one-to-one correspondence with the
$K$-isomorphism classes of objects in ${\mathbf V}$.

We consider the existence of a fine moduli space  as the ideal possibility.
Unfortunately, it is  rarely achieved.
When there is no fine moduli space, we still can ask
for a scheme that best approximates its properties.
Therefore, we  look for schemes $M$ for which there is  a natural
transformation of functors
$$
T_M:Varieties_g(*) \longrightarrow \mor(*, M).
$$
Such schemes certainly exist, for instance, if we work over a field $k$
then  $M=\spec k$. All schemes $M$ for which $T_M$
exists form an inverse system which is closed under fiber products.
Thus,  as long as we are not unlucky, there is a universal (or largest)
scheme with this property. Though it is not usually done, it
should be called the {\it categorical moduli space}.

This object can be rather useless in general.
For instance, fix $n, d$ and let
${\mathbf H}_{n,d}$ be the class of all hypersurfaces of degree $d$
in $\p^{n+1}_k$ up to isomorphisms. One can see
that a categorical moduli space exists and it is $\spec k$.

In order to get something reasonable, we impose extra conditions.
A  scheme ${\rm Moduli}_{\mathbf V}$
is a {\it coarse moduli space} for ${\mathbf V}$ if
the following hold.
\begin{enumerate}
\item There is a natural transformation of functors
$$
{\rm ModMap}:Varieties_{\mathbf V}(*) \longrightarrow
 \mor(*, {\rm Moduli}_{\mathbf V}),
$$
\item ${\rm Moduli}_{\mathbf V}$ is universal satisfying (1), and
\item for any algebraically closed field $K\supset k$,
$$
{\rm ModMap}:Varieties_{\mathbf V}(\spec K) \stackrel{\cong}{\longrightarrow}
 \mor(\spec K, {\rm Moduli}_{\mathbf V})={\rm Moduli}_{\mathbf V}(K)
$$
is an isomorphism (of sets).
\end{enumerate}
 In many cases, the naturally occurring
moduli spaces have a further very useful property.
\begin{enumerate}\setcounter{enumi}{3}
\item There is a $V_U\to U$ in ${\mathbf V}$
such that the corresponding moduli map
$U\to  {\rm Moduli}_{\mathbf V}$ is open and quasi finite.
\end{enumerate}
Following woodworking terminology, I propose to call
a moduli space satisfying conditions (1--4) a
{\it bastard} %or {\it millbastard}   
moduli space.
\end{say}

\begin{say}[Problems with the moduli of smooth varieties]
\label{higher.nongen.mod.say.intro}{\ }

In contrast with curves, the
moduli theory for higher dimensional smooth varieties can be
 very badly behaved, as shown by the following examples.

(\ref{higher.nongen.mod.say.intro}.1)  (Ruled surfaces)
Let $C$ be a smooth curve and  $L$ a line bundle on $C$ that is
generated by 2 sections $f,g$. On $S:=C\times \a^1_t$,
with first projecion $\pi_1$,
consider the exact sequence
$$
0\to \pi_1^*L^{-1}\stackrel{(t,f,g)}{\longrightarrow}
\pi_1^*L^{-1}+\o_S+\o_S\to Q\to 0.
$$
$Q$ is a rank 2 vector bundle on  $C\times \a^1_t$.
We can view $\p_{C\times \a^1_t}Q$ is a $\p^1$-bundle over
$S$, or as a flat family of ruled surfaces over $\a^1$.

If $t\neq 0$ then $t:\pi_1^*L^{-1}\to \pi_1^*L^{-1}$
is an isomorphism, thus $Q_t\cong \o_C+\o_C$.
If $t=0$ then we get
$$
Q_0\cong L^{-1}+\coker\bigl[L^{-1}\stackrel{(f,g)}{\longrightarrow}
\o_C+\o_C\bigr]
\cong L^{-1}+L.
$$ 

Thus we get a flat family of smooth ruled surfaces whose
general member is $\p^1\times C$ and whose special member
is $\p_C(L^{-1}+L)$. 
In a coarse moduli space over $\c$ both of these
should correspond to $\c$-points, but 
the above family shows that  the moduli point
$\bigl[\p_C(L^{-1}+L)\bigr]$ is in the closure of the moduli point
 $\bigl[\p^1\times C\bigr]$. This is impossible 
(at least for schemes but not  for stacks). 

One can be even more specific for $C=\p^1$.
Minimal ruled surfaces over $\p^1$ are
$\f_m:=\p_{\p^1}\bigl(\o_{\p^1}+\o_{\p^1}(m)\bigr)$ for $m\geq 0$. 
The ``moduli space'' has 2 connected components, corresponding to
even and odd values of $m$.

There are no closed points in this ``moduli space;''
the closure of $\bigl[\f_m\bigr]$ consists of all
the points  
$$
\bigl\{\bigl[\f_m\bigr], \bigl[\f_{m+2}\bigr], \bigl[\f_{m+4}\bigr], \dots\bigr\}.
$$

(\ref{higher.nongen.mod.say.intro}.2) (Abelian, elliptic and K3 surfaces)

A general problem in all these cases
is that a typical deformation of such an algebraic surface over $\c$ is a
non-algebraic complex analytic surface. Thus any algebraic theory
captures only a small part of the full  moduli theory.

(\ref{higher.nongen.mod.say.intro}.3)  For Abelian varieties  and for
 K3 surfaces, the moduli spaces look very strange topologically.
For instance,
the 3-dimensional space of Kummer surfaces
is dense in the 20-dimensional space of all K3 surfaces
 \cite{shaf-ps-MR0284440}.

This can be corrected by fixing  a  basis in
$H^2(*, \z)$, but it is not clear how similar tricks work in
general. Also, as it happens already for stable curves, 
we would like to consider families where not all fibers are
homeomorphic to each other. Then it is no clear what one means
by `` fixing a basis in
$H^2(*, \z)$.''

(\ref{higher.nongen.mod.say.intro}.4) (Repeated blow-ups lead to
non-separatedness)

Let $f:X\to B$ be a smooth family of projective surfaces
over a smooth (affine) pointed curve $b\in B$. 
Let $C_1,C_2, C_3\subset X$ be three sections of $f$,
all passing through a point $x_b\in X_b$ 
that intersect pairwise transversally at $x_b$ and are disjoint elsewhere.

Set $X^1:=B_{C_1}B_{C_2}B_{C_3}X$, where we first blow-up
$C_3\subset X$,  then the birational transform of
$C_2$ in $B_{C_3}X$ and finally the birational transform of
$C_1$ in $B_{C_2}B_{C_3}X$. Similarly, set
$X^2:=B_{C_1}B_{C_3}B_{C_2}X$. Since the $C_i$ are
sections, all these blow-ups are
smooth families of projective surfaces over $B$.
It is easy to check  that
\begin{enumerate}
\item all the fibers are smooth, projective surfaces of general type,
\item $X^1\to B$ and $X^2\to B$ are isomorphic over $B\setminus\{b\}$,
\item the fibers $X^1_b$ and $X^2_b$  are  {\em not} isomorphic.
\end{enumerate}

This type of behavior happens
every time we look at deformations of a surface with at least
3 points blown-up.

(\ref{higher.nongen.mod.say.intro}.5) (Non-separatedness
for minimal resolutions.)

Let  $X_0:=\bigl(f(x_1,\dots, x_4)=0\bigr)\subset \p^3$
be a surface of degree $n$ that
has an ordinary double point  at  $p=(0{:}0{:}0{:}1)$
as its sole singularity 
and contains the pair of lines $(x_1x_2=x_3=0)$.
Let $g$ be homogeneous of degree $n-1$ such that  $x_4^{n-1}$ appears in it
with nonzero coefficient.
Consider the family of surfaces
$$
X:=\bigl(f(x_1,\dots, x_4)+tx_3g(x_1,\dots, x_4)=0\bigr)
\subset \p^3_{\mathbf x}\times \a^1_t.
$$
Note that  $X_t$ is smooth for general $t\neq 0$
and $X$ contains the pair of smooth surfaces $(x_1x_2=x_3=0)$.

For $i=1,2$, let $X^i:=B_{(x_i,x_3)}X$ denote the blow-up of
 $(x_i=x_3=0)$ with induced morphisms $\pi_i:X^i\to X$ and
$f_i:X^i\to \a^1$.
There is a natural birational map
$\phi:=\pi_2^{-1}\circ \pi_1:X^1\map X^2$.
Let $B_{p}X$ denote the blow-up of
$p=\bigl((0{:}0{:}0{:}1),0\bigr)$ with exceptional divisor
$E\subset B_{p}X$.
One checks  that
\begin{enumerate}
\item all the fibers are smooth, projective minimal models,
\item $X^1\to B$ and $X^2\to B$ are isomorphic over $B\setminus\{b\}$,
\item the fibers $X^1_b$ and $X^2_b$  are isomorphic, but
\item $X^1\to B$ and $X^2\to B$ are {\em not} isomorphic.
\end{enumerate}

While it is not clear from our construction,
similar problems happen for any smooth family of surfaces
where the general fiber has ample canonical class and a special fiber
has nef (but not ample)  canonical class, see \cite{art-bri, briesk, r-c3f}.

\end{say}

\begin{say}[Answers to these problems]\label{ans.4}
The problems (\ref{higher.nongen.mod.say.intro}.1--3)
 come from the global geometry of the varieties
that we work with. 

The current assumption is that the moduli
problem of uniruled varieties is usually pathological.
Furthermore, any general attempt to create a good moduli functor ends up
with a theory that is not compatible with birational equivalence.
(Although there are examples, like the moduli of
smooth hypersurfaces of degree $n$ in $\p^n$ for $n\geq 4$,
where the biregular moduli theory ends up being birationally invariant.
However, even in these cases, 
it is not clear that a sensible compactification exists.)

For varieties with trivial canonical class one should get a nice moduli theory
only after some suitable ``rigidification.'' This can consist of
choosing a basis in $H^2(*,\z)$ or fixing an ample divisor.
The compactification question is mostly still unsolved.
For instance, a geometrically meaningful compactification of the moduli of
K3 surfaces is yet to be found.

The problems (\ref{higher.nongen.mod.say.intro}.4--5)
are more local. The aim of these notes is to explain
how to deal with them for varieties of general type.
The solution is to work with
{\it canonical models} 
instead of smooth  varieties of general type.

Following
\cite{iitaka}
(see \cite[Sec.2.1.C]{laz-book} for  details) 
a  smooth projective variety $X$ of
dimension $n$ is  of {\it general type}
if  the
following  equivalent conditions hold:

\begin{enumerate}
\item $h^0\bigl(X, \o_X(mK_X)\bigr)\geq \epsilon\cdot m^n$ for
some $\epsilon>0$ and $m\gg 1$.
\item $\proj R(X,K_X)$ has dimension $n$.
\item The natural map $X\map \proj R(X,K_X)$ is birational.
\end{enumerate}

The main reason, however, why we do not study  the 
moduli functor of smooth varieties up to isomorphism is that, 
in dimensions two and up, 
smooth projective  varieties do not form the
{\em smallest} basic class. Given any smooth projective variety
$X$, one can blow up any set of points or subvarieties $Z\subset X$
to get another  smooth projective variety $B_ZX$ which is
very similar to $X$. 
Therefore, the  basic object should be not a single smooth projective   variety
but a whole {\it birational equivalence class} of
smooth projective  varieties. 
Thus it would be better to work with smooth, proper  families $X\to S$
modulo birational equivalence over $S$. 
That is, with the 
moduli functor
$$
GenType_{bir}(S):=
\left\{
\begin{array}{c}
\mbox{Smooth, proper  families $X\to S$,}\\
\mbox{every fiber is of general type,}\\
\mbox{modulo birational equivalence over $S$.}
\end{array}
\right\}
\eqno{(\ref{ans.4}.4)}
$$

In essence this is what we end up doing,
but it is very cumbersome do deal with
birational equivalence over a base scheme.
\end{say}

The following result, or rather, its proof,
shows the way to a good moduli theory for
varieties of general type.

\begin{prop}\label{ample-mm.prop}
 Let $f_i:X^i\to B$ be two smooth families of projective 
varieties over a  smooth curve $B$
such that the canonical classes $K_{X^i}$ are $f_i$-ample.
Then every  isomorphism
between the  generic fibers $\phi:X^1_{k(B)}\cong X^2_{k(B)}$
 extends to an isomorphism $\Phi:X^1\cong X^2$.
\end{prop}

Proof.  Let 
$\Gamma\subset X^1\times_B X^2$ be the closure of the graph of
$\phi$. Let $Y\to \Gamma$ be the normalization, 
with projections $p_i:Y\to X^i$ and  $f:Y\to B$.
We use the canonical class to compare the $X^i$.
Since the $X^i$ are smooth,
$$
K_Y\sim p_i^*K_{X^i}+E_i\qtq{where $E_i$ is effective and  $p_i$-exceptional.}
\eqno{(\ref{ample-mm.prop}.1)}
$$
Since $(p_i)_*\o_Y(mE_i)=\o_{X^i}$ for every $m\geq 0$, we get that
$$
\begin{array}{rclcl}
(f_i)_* \o_{X^i}\bigl(mK_{X^i}\bigr)&=&
(f_i)_* (p_i)_* \o_{X^i}\bigl(mp_i^*K_{X^i}\bigr)&=&\\
&=&
(f_i)_* (p_i)_* \o_{X^i}\bigl(mp_i^*K_{X^i}+mE_i\bigr)&=&\\
&=&
(f_i)_* (p_i)_*\o_Y\bigl(mK_Y\bigr)&=&
f_*\o_Y\bigl(mK_Y\bigr).
\end{array}
$$
Since the $K_{X^i}$ are $f_i$-ample,
$X^i=\proj_B \sum_{m\geq 0} (f_i)_* \o_{X^i}\bigl(mK_{X^i}\bigr)$.
Putting these together, we get the isomorphism
$$
\begin{array}{rclcl}
\Phi: X^1&\cong &\proj_B \sum_{m\geq 0} (f_1)_* \o_{X^1}\bigl(mK_{X^1}\bigr)
&\cong &\\
&\cong &\proj_B \sum_{m\geq 0}f_*\o_Y\bigl(mK_Y\bigr)&\cong &\\
&\cong &
\proj_B \sum_{m\geq 0} (f_2)_* \o_{X^2}\bigl(mK_{X^2}\bigr)&\cong& X^2.
\quad \qed
\end{array}
$$

Note that the smoothness of the $X^i$ is used only through the
pull-back formula (\ref{ample-mm.prop}.1).
This leads to the first major definition:

\begin{prel-defn} \label{prel-defn.defn} Let $B$ be a smooth curve, 
$X$ a normal variety and
$f:X\to B$ a non-constant projective  morphism.
Assume that
\begin{enumerate}
\item $m_0K_X$ is Cartier for some $m_0>0$.
(This is needed to make sense of (2) and also to define the pull-back in (3).)
\item  $m_0K_X$ is  $f$-ample.
\item  If $p:Y\to X$ is a resolution of singularities then
$$
m_0K_Y\sim p^*\bigl(m_0K_{X}\bigr)+E(m_0)
$$ 
where $E(m_0)$ is effective and  $p$-exceptional.
\end{enumerate}
We do not want to keep carrying the $m_0$ along, thus we switch to
$\q$-divisors and $\q$-linear equivalence and write (3)  as
\begin{enumerate}
\item[(3')]  If $p:Y\to X$ is a resolution of singularities then
$$
K_Y\simq p^*K_{X}+E
$$ 
where $E$ is an effective,  $p$-exceptional $\q$-divisor.
\end{enumerate}
With these assumptions we can make the following informal definitions:
\begin{enumerate}\setcounter{enumi}{3}
\item A ``general'' fiber of $f$ is a {\it canonical model.}
\item A ``special'' fiber of $f$ is a {\it semi log canonical model}
if (1--3') continue to hold for every base change $X':=X\times_BB'\to B'$,
for every smooth curve $B'$.
\end{enumerate}
In this area, ``semi'' refers to allowing non-normal schemes
and ``log'' refers to allowing exceptional divisors with coefficients
$\geq -1$, see (\ref{sl-can.def.1sr.intro}).

(Note that (1--2) are inherited by  every base change, thus  the only
question is (3'). Already for curves, this condition is necessary.
Indeed, let $(f(x,y,z)=0)$ be any plane curve with isolated singularities
and $(g(x,y,z)=0)$ a general smooth plane curve. Then
$$
\bigl( f(x,y,z)+tg(x,y,z)=0\bigr)\subset \p^2\times \a^1_t
$$
is smooth near the $t=0$ fiber. After the base change $t=s^m$ 
we get
$$
\bigl( f(x,y,z)+s^mg(x,y,z)=0\bigr)\subset \p^2\times \a^1_s
$$
and now (3') is satisfied for $m\geq 3$ iff $(f(x,y,z)=0)$
has only ordinary nodes.)
\end{prel-defn}

\subsection*{Basic principles}{\ }

The moduli theory of higher dimensional varieties of general type
is governed by the following 
four basic principles.

\begin{principle} \label{princ1}
Canonical models are the correct higher dimensional analogs
of smooth projective curves of genus $\geq 2$.
\end{principle}

\begin{principle} \label{princ2}
Semi log canonical models are the correct higher dimensional analogs
of stable projective curves of genus $\geq 2$.
\end{principle}

\begin{principle} \label{princ3}
Flat families of canonical models form the correct higher dimensional 
open moduli problem.
\end{principle}

\begin{principle} \label{princ4}
Flat families of semi log canonical models do {\bf not}  form 
the correct higher dimensional 
compactified moduli problem.
\end{principle}

\section{Canonical models}\label{sec.2}

As noted in (\ref{princ1}), canonical models are the
basic objects of our moduli theory. 

\begin{defn}  \label{can.def.1sr.intro}
A normal projective  variety $X$ is a {\it canonical model}
 iff the following hold.
\begin{enumerate}
\item $m_0K_X$ is Cartier for some $m_0>0$.
\item for every (equivalently, for one) resolution of singularities $f:X'\to X$ 
there is an effective, $f$-exceptional $\q$-divisor $E$ such that
$$
K_{X'}\simq f^*K_X+E.
$$
\item $K_X$ is ample.
\end{enumerate}
Singularities satisfying (1--2) are called {\it canonical.}
\end{defn}

This assumption implies that the {\it canonical ring} of $X$
$$
R(X,K_X):=\sum_{m\geq 0} H^0\bigl(X, \o_X(mK_X)\bigr)
$$
is isomorphic to the  canonical ring of $X'$
$$
R(X',K_{X'}):=\sum_{m\geq 0} H^0\bigl(X', \o_{X'}(mK_{X'})\bigr).
$$
In particular, 
$$
X=\proj_k R(X,K_X)=\proj_k R(X',K_{X'})
$$
is the unique canonical model in  the birational equivalence class of $X$.

Now we know \cite{bchm, siu-fg} that the
 canonical ring
$R(X,K_{X})$ of a smooth projective variety of general type
is always finitely generated, thus
$X^{can}$ is a projective variety. 
It is not obvious, but true,  that 
$X^{can}$  is a canonical model \cite{r-c3f}.

\begin{defn}[Moduli  functor of canonical models]\label{can.mod.funct.def.intro}
The 
moduli functor of canonical models is
$$
CanMod(S):=
\left\{
\begin{array}{c}
\mbox{Flat, proper  families $X\to S$,}\\
\mbox{every fiber is a canonical model,}\\
\mbox{modulo isomorphisms over $S$.}
\end{array}
\right\}
\eqno{(\ref{can.mod.funct.def.intro}.1)}
$$
This is an improved version of
 the 
birational moduli functor 
$GenType_{bir}(*)$ (\ref{ans.4}.1).

(Traditionally this was considered to be the obviously correct definition, but,
in view of Principle \ref{princ4}, it needs an explanation.
For details, see (\ref{why.can.better}).)

By a theorem of \cite{siu-plur}, in a
smooth, proper  family of varieties of  general type
the canonical rings  form a  flat family
and so do the canonical models. Thus there is a natural transformation
$$
T_{\rm CanMod}: GenType_{bir}(*)\to CanMod(*).
$$
By definition and by (\ref{ample-mm.prop}), if $X_i\to S$ are two 
smooth, proper  families of varieties of  general type
then
$$
T_{\rm CanMod}(X_1/S)=T_{\rm CanMod}(X_2/S)\qtq{iff \quad 
$X_1$ and $X_2$ are birational,}
$$
thus $T_{\rm CanMod}$ is injective. It is, however, not surjective,
but we have the following 
 partial surjectivity statement.

 Let $Y\to S$ be a flat family of canonical models.
Then there is a dense open subset $S^0\subset S$ and a
smooth, proper  family of varieties of  general type $Y^0\to S^0$
such that 
$$
T_{\rm CanMod}(Y^0/S^0)=[X^0/S^0].
$$
Some of the obstruction to surjectivity are obvious but some,
as in (\ref{higher.nongen.mod.say.intro}.3),  are quite subtle.
\end{defn}

Canonical curves are exactly the smooth, projective curves of genus
$\geq 2$. Canonical surfaces  have at most Du~Val singularities
(also called rational double points).
Starting with dimension 3, we get more complicated singularities.
For instance,
$$
\bigl(x_0^{a_0}+\cdots+x_n^{a_n}=0\bigr)\subset \a^{n+1}
$$ is canonical iff $\sum \frac 1{a_i}>1$
and a quotient of $\a^n$ be a finite subgroup $G$ 
without quasi-reflections is canonical iff
for every $g\in G$ ($g\neq 1$)
with eigenvalues $e^{2\pi i c_j}$ (where $0\leq c_j<1$) we have
$\sum_j c_j\geq 1$.

The most important general property of canonical  singularities
is the following. For  short proofs see  \cite[Sec.5.1]{km-book}
or \cite{k-kv}.

\begin{thm}\label{klt.rtl.thm}\cite{elkik}
Let $X$ be a canonical model  over a field of characteristic 0.
Then $X$ has  rational singularities.
That is, 
$R^if_*\o_Y=0$ for $i>0$ for every resolution of singularities
 $f:Y\to X$.
\end{thm}

Using the covering trick of \cite{r-c3f} this implies that
the reflexive hulls
$\omega_X^{[m]}$ are CM (=Cohen-Macaulay) for every $m$.

\section{Semi log canonical models}\label{sec.3}

First we translate (\ref{prel-defn.defn}.5)
into a proper definition.

 Let $B$ be a smooth curve, 
$X$ a normal variety and
$f:X\to B$ a non-constant projective  morphism.
When is the fiber $X_b$ a semi log canonical model?

By \cite{kkms} there is a smooth pointed curve $b'\in B'$ and a
finite morphism $(b'\in B')\to (b\in B)$ such that the base change
$f':X':=X\times_BB'\to B'$ has a resolution
$\pi:Y'\to X'$ such that
$(f'\circ \pi)^{-1}(b')\subset Y'$ is a reduced
simple normal crossing divisor.
We can also assume that the birational transform
$Y'_{b'}:=\pi^{-1}_*X'_{b'}$ is smooth.
Thus $
(f'\circ \pi)^{-1}(b')= Y'_{b'}+F$ where $F$ is a reduced
simple normal crossing divisor.
By (\ref{prel-defn.defn}.3), 
$$
K_{Y'}\simq {\pi}^*K_{X'}+E'
$$ where $E'$ is effective.
%Any fiber of $f'\circ \pi$ is linearly equivalent to 0, hence 
Since $Y'_{b'}+F'=\pi^*\bigl(X'_{b'}\bigr)$, we can
rewrite the above as
$$
K_{Y'}+Y'_{b'}\simq {\pi}^*\bigl(K_{X'}+X'_{b'}\bigr)+E'-F.
$$
Restricting to $Y'_{b'}$ we get
$$
 K_{Y'_{b'}}\simq 
\pi^*K_{X'_{b'}}+\bigl(E'-F\bigr)|_{Y'_{b'}}.
$$
Since $E'$ is effective, its restriction is again effective,
but the restriction of $-F$ brings in
negative coefficients. However,  none of these is smaller than $-1$
since $F$ and  $Y'_{b'}$ intersect transversally.

\begin{defn}  \label{sl-can.def.1sr.intro}
A reduced, projective  variety $X$ is a {\it semi log canonical model} 
or {\it slc model} iff
the following hold.
\begin{enumerate}
\item $m_0K_X$ is Cartier  for some $m_0>0$.
\item $X$ has only ordinary nodes in codimension 1.
\item For every resolution of singularities $f:X'\to X$ 
there is an  $f$-exceptional $\q$-divisor $E=\sum_i a_iE_i$ such that
$$
K_{X'}\simq f^*K_X+\tsum_i a_iE_i
\qtq{and $a_i\geq -1$ for every $i$.}
$$
\item $K_X$ is  ample. 
\end{enumerate}
One should think of this as combining a global condition (4)
with  purely local conditions (1--3). Singularities satisfying
(1--3) are called {\it  semi log canonical} or {\it slc}.

For slc models it is usually better to use
semi resolutions, that is, 
a proper birational morphism $g:X^s\to X$ such that
$X^s$ has only double normal crossing points $(xy=0)\subset \c^{n+1}$
and pinch points $(x^2=y^2z)\subset \c^{n+1}$ and
$g$ maps the double locus of $X^s$ birationally on the
double locus of $X$. Let  $E$ denote the (reduced) exceptional divisor of $g$.
Then the  canonical ring of $X$
$$
R(X,K_X):=\sum_{m\geq 0} H^0\bigl(X, \o_X(mK_X)\bigr)
$$
is isomorphic to the  {\it semi log canonical ring} of $X^s$
$$
R(X^s,K_{X^s}+E):=\sum_{m\geq 0} H^0\bigl(X^s, \o_{X^s}(mK_{X^s}+mE)\bigr).
$$
This actually creates a lot of problems since 
semi log canonical rings are not always finitely generated
\cite{k-ncsurf}.
\end{defn}

It is a quite subtle theorem that 
semi log canonical models actually satisfy the
preliminary definion (\ref{prel-defn.defn}.5).
This is proved in \cite[17.4]{k-etal} and 
\cite{Kawakita}.
\medskip

To get a feeling for  semi log canonical,
let us review the classification of slc surface singularities.

\medskip
\subsection*{Singularities of semi log canonical  surfaces}{\ }
\medskip

It is convenient to describe the singularities of  log canonical  surfaces
by the dual graph of their minimal resolution.
That is, given a singularity $(s\in S)$
with minimal resolution $g:X\to S$
we draw a graph  $\Gamma$ whose vertices are the $g$-exceptional curves
and two vertices are connected by an edge iff the corresponding curves
intersect. 
We use the number $-(E_i\cdot E_i)$ to represent a vertex.
In our examples, save in (\ref{lcnltnb}.4.a), 
 all the exceptional curves are isomorphic to $\p^1$.

Let $\det(\Gamma)$ denote the determinant of the negative of the
intersection matrix of the dual graph. This matrix
is 
positive definite for exceptional curves. 
For instance, if $\Gamma=\{ 2 \ - \ 2  \ - \ 2 \}$
then
$$
\det(\Gamma)=
\det\left(
\begin{array}{rrr}
2 & -1 & 0\\
-1 & 2 & -1\\
0 & -1 & 2
\end{array}
\right)=4.
$$

\begin{say}[List of log canonical surface singularities]{\ }
\label{lcnltnb}

Each case includes all previous ones.

(\ref{lcnltnb}.1) Terminal = smooth.

(\ref{lcnltnb}.2) Canonical = Du~Val (= rational double point).

(\ref{lcnltnb}.3) Log terminal = quotient of $\c^2$ by a finite subgroup
of $GL(2,\c)$ that acts freely outside the origin.
The order of the group is $\det(\Gamma)$.
A more detailed list is the following:

(a)  (Cyclic quotient)
$$
 c_1  \ -\ \cdots  \ - \ c_n 
$$ 

(b) (Dihedral quotient) Here $n\geq 2$ with dual graph
$$
\begin{array}{ccccccc}
 &&&&&& 2\\
&&&&& \diagup &\\
  c_1 & - & \cdots  & - & c_n &&\\
&&&&& \diagdown &\\
 &&&&&& 2
\end{array}
$$

(c) (Other quotients) The dual graph has 1 fork
$$
\begin{array}{ccccc}
\Gamma_1 & - & c_0 & - & \Gamma_2\\
 && \vert && \\
&& \Gamma_3
\end{array}
$$
whith 3  cases for 
$\bigl(\det(\Gamma_1),\det(\Gamma_2),\det(\Gamma_3)\bigr)$:
\begin{enumerate}
\item[](Tetrahedral)  (2,3,3)
\item[](Octahedral) (2,3,4)
\item[](Icosahedral)  (2,3,5).
\end{enumerate}

(\ref{lcnltnb}.4) Log canonical 

(a) (Simple elliptic)  $\Gamma=\{E\}$ has a single vertex which
is a smooth elliptic curve
with self intersection $\leq -1$. 

(b) (Cusp) $\Gamma$ is  a circle of smooth rational curves, 
 at least
one of them with with $c_i\geq 3$. 
(The cases $n=1,2$ are somewhat special.)
$$
\begin{array}{ccccccccc}
&& c_n & - &\cdots & - & c_{r+1} && \\
 & \diagup &&&&&& \diagdown &  \\
c_1 &&&&&&&& c_r\\
& \diagdown &&&&&& \diagup &\\
 && c_2 & - &\cdots & - & c_{r-1}&& 
\end{array}
$$

 (c) ($\z/2$-quotient of a cusp or simple elliptic)  
$\Gamma$ has 2 forks.
$$
\begin{array}{ccccccccc}
2 &&&&&&&& 2\\
& \diagdown &&&&&& \diagup &\\
&& c_1 & - & \cdots  & - & c_n &&\\
& \diagup &&&&&& \diagdown &\\
2 &&&&&&&& 2
\end{array}
$$

(d) (Other quotients of a simple elliptic)  
The dual graph is  as in (\ref{lcnltnb}.3.c)
with 3  possibilities for 
$\bigl(\det(\Gamma_1),\det(\Gamma_2),\det(\Gamma_3)\bigr)$:
\begin{enumerate}
\item[] ($\z/3$-quotient)  (3,3,3)
\item[] ($\z/4$-quotient)  (2,4,4)
\item[] ($\z/6$-quotient)   (2,3,6).
\end{enumerate}
\end{say}

If $X$ is a non-normal semi log canonical surface singularity,
then we describe its normalization $\bar X$ together with
the preimage of the double curve $\bar B\subset \bar X$.

The {\it extended dual graph}  $(\Gamma, \bar B)$ has an additional vertex
(repesented by  $\bullet$) for each local branch of $\bar B$ 
connected to $C_i$ if $(\bar B\cdot C_i)\neq 0$.

\begin{say}[List of semi log canonical surface singularities]
\label{surf.with.B.list} There are 3 irreducible cases.
(The number on some edges is the different, which we do not
define here \cite[Sec.16]{k-etal}. Their
role is explained in (\ref{surf.with.B.list}.4).

(\ref{surf.with.B.list}.1)  (Cyclic quotient, one branch of $\bar B$)
$$
\bullet   \ \stackrel{1-\frac1{\det\Gamma}}{-\!\!\!-\!\!\!-\!\!\!-} 
\ c_1  \ -\ \cdots  \ - \ c_n  
$$ 

(\ref{surf.with.B.list}.2)  (Cyclic quotient, two branches of $\bar B$)
$$
\bullet   \ \stackrel{1}{-}\ c_1  \ -\ \cdots  \ - \ c_n  \ \stackrel{1}{-}
\ \bullet 
$$ 

(\ref{surf.with.B.list}.3) (Dihedral quotient) Here $n\geq 2$ with dual graph
$$
\begin{array}{ccccccc}
 &&&&&& 2\\
&&&&& \diagup &\\
\bullet   \ \stackrel{1}{-}\  c_1 & - & \cdots  & - & c_n &&\\
&&&&& \diagdown &\\
 &&&&&& 2
\end{array}
$$

(\ref{surf.with.B.list}.4) (Reducible cases)
We can take several components as above and glue them together
along two local branches of $\bar B$. The gluing is allowed
only if we see the same numbers on the edges.

Thus we can glue 2 copies as in (\ref{surf.with.B.list}.1)
 as long as both have the same $\det(\Gamma)$.

Or we can take any number of those in (\ref{surf.with.B.list}.2),
make a chain out of them and then either turn the chain
into a circle or end it with a copy of (\ref{surf.with.B.list}.3).

We are also allowed to glue a local branch of $\bar B$
to itself by an involution. 
For instance, $\bullet   \ - \ 1$ glued to itself
gives the pinch point $(x^2=y^2z)\subset \a^3$. 
 
\end{say}

Note: The above dual graphs are correct in any characteristic,
the descriptions as quotients are correct as long as
the the characteristic does not divide the order of the
group mentioned.

\medskip
\subsection*{Du~Bois singularities}{\ }
\medskip

Semi log canonical singularities need not be rational,  not even CM
(=Cohen-Macaulay) and their most important property is that they are
Du~Bois. After some examples, we discuss Du~Bois singularities
and their useful properties.

\begin{exmp}
It is easy to see that a cone over a smooth variety $X\subset \p^N$ 
is log canonical iff $K_X\simq r\cdot H$ for some $r\leq 0$
where $H$ is the hyperplane class. For us the interesting case is
when $K_X\simq 0$ (hence $r=0$).
For these, the cone is CM (resp. rational)
iff $H^i(X, \o_X)=0$ for $0<i<\dim X$
(resp.\ for $0<i\leq \dim X$). Thus we see the following.
\begin{enumerate}
\item A cone over an Abelian variety $A$ is CM iff $\dim A=1$.
\item  A cone over a K3 surface is CM but not rational.
\item  A cone over an Enriques surface is CM and rational.
\item  A cone over a smooth Calabi-Yau complete intersection is CM 
but not rational.
\end{enumerate}
\end{exmp}

The concept of Du~Bois singularities
was introduced by
Steenbrink in \cite{Steenbrink83} as a weakening of rationality. 
The precise definition is rather
involved, but our main applications rely only on the following
consequence.

\begin{thm}\label{lc.is.db.cor}\cite{k-db}, \cite[Chap3]{k-m-book}
  Let $X$ be a proper slc scheme  over $\c$.
 Then the natural map
  $$
  H^i(X^{\rm an},\c)\to H^i(X^{\rm an},\o_{X^{\rm an}})\cong  H^i(X,\o_{X})
  $$
  is surjective for all $i$.
(In fact, with a functorial splitting.)
\end{thm}

In studying moduli questions, it is very useful to know
that certain numerical invariants are locally constant.
All of these follow from the  Du~Bois property,
via the following base-change theorem 
\cite{dub-jar, dub}.

\begin{prop}\label{DB-J.coh.inv}%\cite{MR0376678}
  Let $f:X\to S$ be a flat, proper morphism.
Assume that the fiber $X_s$ 
is Du~Bois  for some $s\in S$.
 Then there is an open  neighborhood  $s\in S^0\subset S$
such that, for all $i$,
\begin{enumerate}
\item $R^if_*\o_{X}$ is locally
  free  and compatible with  base change 
over  $S^0$ and
\item  $s\mapsto  H^i(X_s, \o_{X_s}\bigr)$ is a 
locally constant function on  $S^0$.
\end{enumerate}
\end{prop}

Proof.  By Cohomology and Base Change
\cite[III.12.11]{hartsh}, the theorem is
equivalent to proving that the restriction maps
$$
\phi^i_s:R^if_*\o_{X}\to H^i(X_s, \o_{X_s}\bigr)
\eqno{(\ref{DB-J.coh.inv}.3)}
$$ are surjective for every $i$. By the Theorem on Formal Functions
\cite[III.11.1]{hartsh}, it is enough to prove this when
$S$ is replaced by any $0$-dimensional scheme $S_n$ whose closed point is
$s$.

Thus assume from now on that we have  a flat, proper morphism $f_n:X_n\to S_n$,
$s\in S_n$ is the only closed point and $X_s$ is Du~Bois.
Then $H^0\bigl(S_n,  R^if_*\o_{X}\bigr)=H^i\bigl(X_n, \o_{X_n}\bigr)$,
hence we can identify the  $\phi^i_s$ with the maps
$$
\psi^i:H^i\bigl(X_n, \o_{X_n}\bigr)\to H^i(X_s, \o_{X_s}\bigr)
\eqno{(\ref{DB-J.coh.inv}.4)}
$$

By the Lefschetz principle  we may assume that
everything is defined over $\c$.
By GAGA  (cf.\ \cite[App.B]{hartsh}), 
both sides of (\ref{DB-J.coh.inv}.4) are unchanged if
we replace $X_n$ by the corresponding analytic space $X^{\rm an}_n$.
Let $\c_{X_n}$ (resp.\ $\c_{X_s}$) denote the sheaf of locally constant
fucntions on $X_n$ (resp.\ $X_s$) and
$j_n:\c_{X_n}\to \o_{X_n}$ (resp.\  $j_s:\c_{X_s}\to \o_{X_s}$)
the natural inclusions. We have a commutative diagram
$$
\begin{array}{ccc}
H^i\bigl(X_n, \c_{X_n}\bigr)&\stackrel{\alpha^i}{\to} 
&H^i(X_s, \c_{X_s}\bigr)\\[1ex]
j^i_n\downarrow\hphantom{j^i_n} && \hphantom{j^i_s}\downarrow j^i_s\\
H^i\bigl(X_n, \o_{X_n}\bigr)&\stackrel{\psi^i}{\to}& H^i(X_s, \o_{X_s}\bigr)
\end{array}
$$
Note that $\alpha^i$ is an isomorphism since
the inclusion $X_s\into X_n$ is a homeomorphism and
$j^i_s$ is surjective since $X_s$ is Du~Bois.
Thus $\psi^i$ is also surjective.\qed
\medskip

A  line bundle $L$ on $X$ is called  {\it $f$-semi ample}
if there is an $m>0$ such that $L^m$ is $f$-generated by global sections.
Using cyclic coverings \cite[Sec.2.4]{km-book}, there is a
finite
morphism $\pi:Y\to X$ such that
$\pi_*\o_Y=\sum_{r=0}^{m-1} L^{-r}$ and
$f\circ\pi: Y\to S$
 also  has  slc fiber over $s$.
Thus (\ref{DB-J.coh.inv}) implies the following.

\begin{cor}
  \label{semiamp.coh-inv.cor}
Let $f:X\to S$ be a proper and  flat morphism  with 
 slc fibers  over closed points; $S$
  connected.   
Let $L$ be an  $f$-semi ample line bundle on $X$.
Then,  for all $i$,
\begin{enumerate}
\item $R^if_*\bigl(L^{-1}\bigr)$ is locally
  free and compatible with base change and 
\item  $H^i(X_s, L_{s}^{-1}\bigr)$ is independent of $s\in S$. \qed
\end{enumerate}
\end{cor}

Choose $L$ to be $f$-ample above. 
By \cite[5.72]{km-book},
  $X_s$ is CM iff $H^i(X_s, L_{s}^{-m}\bigr)=0$ for all $m\gg 1$
and $i<\dim X$. The latter properties are deformation invariant
for slc fibers by (\ref{semiamp.coh-inv.cor}). Thus we conclude:

\begin{cor}\label{cm.def.inv}
  \label{cm-inv.cor} 
 Let $f:X\to S$ be a projective and  flat morphism  with 
 slc fibers   over closed points; $S$
  connected.
  Then, if one fiber of $f$ is CM  then
  all fibers of $f$ are CM. \qed
\end{cor}

(Nore that for arbitrary flat, projective morphisms $f:X\to S$, the set of 
points $s\in S$ such that the fiber $X_s$ is CM 
is open, but usually not closed.) 

The next example shows that
non-CM varieties occur among the
irreducible components of smoothable, CM and slc varieties.

\begin{exmp}\label{noncm.6.exmp}
Here is an example of 
a stable
family of  projective varieties $\{Y_t: t\in T\}$
such that
\begin{enumerate}
\item $Y_t$ is smooth, projective for $t\neq 0$, 
\item   $K_{Y_t}$  is ample and Cartier for every $t$,
\item $Y_0$ is slc and CM,
\item the irreducible components of $Y_0$ are
normal, but 
\item one of the irreducible components of $Y_0$ is not CM.
\end{enumerate}

Let $Z$ be a smooth Fano variety 
of dimension $n\geq 2$ such that
$-K_Z$ is very ample, for instance $Z=\p^2$.
Set $X:=\p^1\times Z$ and view it as embedded
by $|-K_X|$ into $\p^N$ for suitable $N$.
Let $C(X)\subset \p^{N+1}$ be the cone over $X$.

Let $M\in |-K_Z|$ be a smooth member and consider
the following divisors in $X$:
$$
D_0:=\{(0:1)\}\times Z,\
D_1:=\{(1:0)\}\times Z\qtq{and}
D_2:=\p^1\times M.
$$
Note that $D_0+D_1+D_2\sim -K_X$.
Let $E_i\subset C(X)$ denote the cone over $D_i$.
Then $E_0+E_1+E_2$
is a hyperplane section of $C(X)$ and 
$\bigl(C(X), E_0+E_1+E_2\bigr)$ is lc. 

For some $m>0$, let $H_m\subset C(X)$ be a 
general intersection with a degree $m$ hypersurface.
Then 
$$
\bigl(C(X), E_0+E_1+E_2+H_m\bigr)
$$
is snc outside the vertex and is lc at the vertex.
Set $Y_0:=E_0+E_1+E_2+H_m$
Since $Y_0\sim \o_{C(X)}(m+1)$, 
we can view it as a slc limit of a
family of smooth hypersurface sections $Y_t\subset C(X)$.

The cone over $X$ is CM, hence
its hyperplane section  $E_0+E_1+E_2+H_m$ is also CM.
However, $E_2$ is not CM.
To see this, note that
$E_2$ is the cone over $\p^1\times M$ and,
by the K\"uneth formula,
$$
H^i(\p^1\times M, \o_{\p^1\times M})=H^i(M, \o_M)=
\left\{
\begin{array}{l}
k\qtq{if $i=0,n-1$,}\\
0 \qtq{otherwise.}
\end{array}
\right.
$$
Thus $E_2$ is not CM.
\end{exmp}

As in the proof of \cite[III.9.9]{hartsh},
we get from (\ref{semiamp.coh-inv.cor}) the following.

\begin{prop}
  \label{omega.exists.pot.slc.prop}
Let $f:X\to S$ be a projective,  flat morphism  with 
 slc fibers  over closed points.
Then  $\omega_{X/S}$ exists and is compatible with base change.
That is,  for  any $g:T\to S$ 
the natural map
$$
g_X^*\omega_{X/S}\to \omega_{X_T/T}
\qtq{is an isomorphism}
$$
where  $g_X:X_T:=X\times_ST\to X$ is the first projection. \qed
\end{prop}

(This seems like a very complicated way to prove that
$\omega_{X/S}$ behaves as expected, but, as far as I can tell, this 
was not known before. A proof for non-projective algebraic maps
is given in \cite{k-kv}. I do not know how to prove 
(\ref{omega.exists.pot.slc.prop}) for analytic morphisms  $f:X\to S$.)

If the fibers $X_s$ are CM, then
$H^i(X_s, \omega_{X_s}\otimes L_s\bigr)$ is dual to
$H^{n-i}(X_s, L_s^{-1}\bigr)$, and the following is clear.
In general, a more detailed inductive argument is needed
\cite[Chap.4]{k-m-book}.

\begin{cor}
  \label{semiamp.coh-inv.omega.cor}
Let $f:X\to S$ be a projective,  flat morphism  with 
 slc fibers  over closed points; $S$
  connected.   
Then,  for all $i$,
\begin{enumerate}
\item $R^if_*\omega_{X/S}$ is locally
  free  and compatible with base change and 
\item  $H^i(X_s, \omega_{X_s}\bigr)$ is independent of $s\in S$.\qed
\end{enumerate}
\end{cor}

\section{Moduli of semi log canonical models}\label{sec.4}

Let us illustrate (\ref{princ4}) with the an  example
of a flat, projective  family of surfaces with log canonical
singularties over the pair of lines $(xy=0)\subset \c^2$
such that over one line we have smooth surfaces with ample
canonical class and over the other line we have smooth elliptic surfaces.

\begin{exmp}[Jump of Kodaira dimension]\label{deg.4.surf.p5.intro}

There are 2 families of
nondegenerate degree 4 smooth surfaces in $\p^5$.

One family consists of   Veronese surfaces
$\p^2\subset \p^5$ embedded by $\o(2)$.
The general member of the other family is
$\p^1\times \p^1\subset \p^5$
embedded by $\o(2,1)$, special members are embeddings of the ruled surface
$\f_2$. The two families are distinct since
$$
K_{\p^2}^2=9\qtq{and} K_{\p^1\times \p^1}^2=8.
$$
For both of these surface, a smooth hyperplane section
gives a degree 4 rational normal curve
in $\p^4$.

Let $T_0\subset \p^5$ be 
the cone over the degree 4 rational normal curve
in $\p^4$. $T_0$ has a log canonical (even log terminal)
singularity and 
 $K_{T_0}^2=9$.

For us the intersting feature is that
one  can write $T_0$ as a limit of smooth surfaces
in two distinct ways, corresponding to the two possibilities 
of writing the degree 4 rational normal curve
in $\p^4$  as a hyperplane section of a surface.

From the first family, we get $T_0$ as the special fiber of a flat
family whose general fiber is $\p^2$.
 This family is denoted by $\{T_t:t\in \c\}$.
From the second family,  we get $T_0$  as the special fiber of a flat
family whose general fiber is
 $\p^1\times \p^1$.
 This family is denoted by $\{T'_t:t\in \c\}$.
(In general, one needs to worry about the possibility of
getting embedded points at the vertex, but 
 in both cases the special fiber is indeed $T_0$.)

Note that $K^2$ is constant in the family $\{T_t:t\in \c\}$
but jumps at $t=0$ in the family $\{T'_t:t\in \c\}$.

Next we take a suitable cyclic cover of the two families
to get similar examples with ample canonical class.

Let $\pi_0:S_0\to T_0$ be a double cover, ramified along a
smooth quartic hypersurface section. 
Note that $K_{T_0}\simq -\tfrac32 H$ where $H$ is the hyperplane class.
Thus, by the Hurwitz formula,
$$
K_{S_0}\simq\pi_0^*\bigl(K_{T_0}+2H\bigr)\simq\tfrac12 \pi_0^*H.
$$ 
So $S_0$ has ample canonical class and 
$K_{S_0}^2=2$. Since $\pi_0$ is \'etale over the vertex of $T_0$,
 $S_0$ has 2 singular points,
locally (in the analytic or \'etale topology) isomorphic to the
singularity on $T_0$. Thus  $S_0$ is a log canonical surface.

Both of the smoothings 
lift to  smoothings of $S_0$.
 
From the family $\{T_t: t\in \c\}$ we get a smoothing $\{S_t:t\in \c\}$ where
$\pi_t:S_t\to \p^2$ is a double cover, ramified along
a smooth octic. Thus $S_t$ is smooth,
$K_{S_t}\simq\pi_t^*\o_{\p^2}(1)$
is ample  and $K_{S_t}^2=2$. 

From the family $\{T'_t:t\in\c\}$ we get a smoothing $\{S'_t:t\in \c\}$ where
$\pi'_t:S'_t\to \p^1\times \p^1$ is a double cover, ramified along
a smooth curve of bidegree $(8,4)$. 
One of the families of lines on $\p^1\times \p^1$
pulls back to an elliptic pencil on 
$S'_t$ and $K_{S'_t}^2=0$. 
\end{exmp}

In order to exclude such examples, we
concentrate on the Hilbert function of a slc model.

\begin{defn}[Hilbert function of slc models]
Let $X$ be an slc model. 
 Note that $\omega_{X}$ is locally free outside a
subscheme $Z\subset X$ such that $Z$ has codimension $\geq 2$.
Hence the  reflexive hull
$\omega_{X}^{[m]}:=\bigl(\omega_{X}^{\otimes m}\bigr)^{**}$
is isomorphic to  $\omega_{X}^{\otimes m}$  
over $X\setminus Z$.
The {\it Hilbert function} of
$X$ is
$$
H_X(m):=\chi\bigl(X, \omega_X^{[m]}\bigr).
$$
If $\omega_X^{[N]}$ is locally free, then
$$
\omega_X^{[m_0+mN]}\cong \omega_X^{[m_0]}\otimes \omega_X^{[mN]},
$$
thus  $H_X(m_0+mN)$ is a polynomial in $m$.
Thus we can view $H_X(m)$ as a polynomial in $m$ whose
coefficients are periodic functions (with period $N$).
\end{defn}

We view $\chi\bigl(X, \omega_X^{[m]}\bigr)$ as the basic numerical
invariants of $X$. It is then natural to insist
that they stay constant in ``good'' families of
slc models. Over a reduced base, this is enough to
get the correct definition.

\begin{defn} [Moduli of slc models  over reduced bases]
\label{can.mod.funct.hilb.def}
 Let  $H(m)$ be an integer valued function.
On {\em reduced} schemes, 
the 
moduli functor of semi log canonical models  with  Hilbert function $H$ is
$$
SlcMod_H(S):=
\left\{
\begin{array}{c}
\mbox{Flat, proper  families $X\to S$, fibers are  slc models with}\\
\mbox{ample canonical class and Hilbert function $H(m)$,}\\
\mbox{modulo isomorphisms over $S$.}
\end{array}
\right\}
%\eqno{(\ref{can.mod.funct.hilb.def}.1)}
$$
\end{defn}

Over an arbitrary base, let $f:X\to S$ be a flat, proper family of
slc models. Note that $\omega_{X/S}$ is locally free outside a
subscheme $Z\subset X$ such that $Z\cap X_s$ has codimension $\geq 2$
in each fiber. Each $\omega_{X/S}^{\otimes m}$ is also locally free
on $X\setminus Z$, hence it has a  reflexive hull
$\omega_{X/S}^{[m]}$.

If $s\in S$ is a general point, then
$ \omega_{X/S}^{[m]}|_{X_s}\cong \omega_{X_s}^{[m]}$
but for an abitrary $s\in S$ we only have a restriction map
$$
r_s: \omega_{X/S}^{[m]}|_{X_s}\to \omega_{X_s}^{[m]}
$$
which is, in general, neither injective nor surjective.
The best way to ensure that every fiber of $X_s$ has the same
Hilbert function is to require these restriction maps to be
isomorphisms for every $s\in S$. 
(It turns out that this is the only way, that is, the kernel and the cokernel
of $r_s$ can nor cancel each other for every $s$, unless they are both zero.)
This leads to our final definition.

\begin{defn} [Moduli of slc models]
\label{can.mod.funct.hilb.2.def}
 Let  $H(m)$ be an integer valued function.
The 
moduli functor of semi log canonical models  with  Hilbert function $H$ is
$$
SlcMod_H(S):=
\left\{
\begin{array}{c}
\mbox{Flat, proper  families $X\to S$, fibers are slc models with}\\
\mbox{ample canonical class and Hilbert function $H(m)$,}\\
\mbox{$\omega_{X/S}^{[m]}$ is flat over $S$ and commutes with base change,}\\
\mbox{modulo isomorphisms over $S$.}
\end{array}
\right\}
%\eqno{(\ref{can.mod.funct.hilb.def}.1)}
$$
\end{defn}

\begin{aside}\label{why.can.better} We can now  explain
Principle \ref{princ3}. The reason is that for
flat families of canonical models, 
$\omega_{X/S}^{[m]}$ is automatically flat over $S$ and commutes with base change.
This follows from two special properties of canonical singularities.
For simplicity, consider a flat family $X\to \spec k[\epsilon]$
whose special fiber $X_0$ is affine.

First we use that, as a result of the classification
of canonical surface singularities (\ref{lcnltnb}.2), there is an open
subset $j:U\into X$ whose complement $Z$
 has  codimension $\geq 3$ such that
$\omega_{U_0}$ is locally free.
Thus we have an exact sequence
$$
0\to \epsilon\cdot \omega_{U_0}^m\to \omega_{U}^m\to \omega_{U_0}^m\to 0.
$$
By pushing it forward, we get
$$
0\to \epsilon\cdot j_*\omega_{U_0}^m\to j_*\omega_{U}^m\to j_*\omega_{U_0}^m\to 
\epsilon\cdot R^1j_*\omega_{U_0}^m
$$
As noted after (\ref{klt.rtl.thm}), 
$\omega_X^{[m]}$ is a CM sheaf, hence has depth $\geq 3$ at every point of $Z$.
Therefore, 
$$
R^1j_*\omega_{U_0}^m=H^1\bigl(U_0, \omega_{U_0}^m\bigr)=
H^2_{Z_0}\bigl(X_0, \omega_{X_0}^{[m]}\bigr)=0.
$$
This implies that $\omega_{X}^{[m]} $ equals 
$j_*\omega_{U}^m$ and it is flat over $k[\epsilon]$.
\end{aside}

Now that we have the correct definition, we need to prove
that the corresponding  deformation theory is
 reasonable. The key result is the following.

\begin{thm} \label{omega.all.comm.cor}
Let $f:X\to S$ be  flat, projective morphism 
whose fibers are slc models. Let  $H$ be an integer valued function.

Then there is a locally closed embedding $S_H\into S$
such that a morphism $g:T\to S$ factors through $S_H$
iff  $X\times_ST\to T$ is in $SlcMod_H(T)$.
\end{thm}

For surfaces, a proof of this is outlined in
\cite{hacking}, a general solution is in
\cite{abr-hass}. The following general theory of hulls 
\cite{k-hh} applies in many similar contexts as well.

\begin{defn}\label{hull.field.defn}
  Let $X$ be a scheme over a field $k$ and $F$ a  coherent 
sheaf on $X$.  Set $n:=\dim\supp F$.  
The  {\it hull}  of $F$ is the    unique 
$q:F\to F^{[**]}$ such that
\begin{enumerate}
\item $q$ is an isomorphism at all generic points of $\supp F$,
\item  $q$ is surjective at all codimension 1  points of $\supp F$,
\item  $F^{[**]}$ is $S_2$.
 \end{enumerate}

If $X$ itself is normal, $F$ is coherent   and $\supp F=X$, then
$F^{[**]}$ is the usual reflexive hull $F^{**}$ of $F$.
The hull of a nonzero sheaf is also nonzero, in contrast with the
reflexive hull which kills all torsion sheaves.
\end{defn}

One can construct  $F^{[**]}$ as follows.
First replace $F$ by $F/\tors_{n-1}(F)$
where $\tors_{n-1}(F)$ is the largest subsheaf whose support has
dimension $\leq n-1$.
 Then there is a closed subscheme
$Z\subset \supp F$ of codimension $\geq 2$ such that $F/\tors(F)$ 
 is $S_2$ on $X\setminus Z$. Let $j:X\setminus Z\to X$ be the open embedding
and take
$$
F^{[**]}=j_*\Bigl(\bigl(F/\tors(F)\bigr)|_{X\setminus Z}\Bigr).
$$

\begin{defn}\label{hull.defn}
  Let $f:X\to S$ be a morphism and $F$ a  coherent 
sheaf. 
A {\it hull} of $F$ is a  coherent sheaf $G$ together with a map  $q:F\to G$
such that, 
\begin{enumerate}
\item $G$ is flat over $S$  and
\item for every $s\in S$,
 the induced map 
$q_s:F_s\to G_s$ is a hull (\ref{hull.field.defn}).
\end{enumerate}
It is easy to see  that a hull is unique if it exists.

It is clear from the definition that hulls are preserved by base change.
That is, if $g:T\to S$ is a morphism, $X_T:=X\times_ST$ and
$g_X:X_T\to X$ the first projection then
$g_X^*q:g_X^*F\to g_X^*G$ is also a hull.
\end{defn}

\begin{defn}  Let $f:X\to S$ be a projective morphism and 
 $F$  a  coherent sheaf on $X$. 
%which is
%generically flat on every fiber of $\supp F\to S$.
For a scheme $g:T\to S$ set
 ${\it Hull}(F)(T)=1$ if $g_X^*F$ has a hull
and ${\it Hull}(F)(T)=\emptyset$ if $g_X^*F$ does not have a hull,
where $g_X:T\times_SX\to X$ is the projection.
\end{defn}

The main existence  theorem  is the following.

\begin{thm}[Flattening decomposition for hulls]\label{hull.exists.thm}
  Let $f:X\to S$ be a projective morphism and 
 $F$  a  coherent sheaf on $X$.
%which is generically flat
% on every fiber of $\supp F\to S$.
Then 
\begin{enumerate}
\item ${\it Hull}(F)$  has a fine moduli space ${\rm Hull}(F)$.
\item  The  structure map 
 $\eta:{\rm Hull}(F)\to S$ is a locally closed decomposition, that is,
$\eta$ is one-to-one and onto on geometric points and a locally closed embedding
 on every connected component.
\end{enumerate}
\end{thm}

Applying (\ref{hull.exists.thm}) to the relative dualizing sheaf
gives the following result.

\begin{cor} \label{omega.m.comm.cor}
Let $f:X\to S$ be  projective and equidimensional.
Assume that there is a closed subscheme  $Z\subset X$ such that
 $\codim (X_s, Z\cap X_s)\geq 2$
 for every $s\in S$,   $(X\setminus Z)\to S$ is flat 
and $\omega_{X/S}$ is locally free on $X\setminus Z$.
Then, for any $m$ there is a locally closed decomposition $S_m\to S$
such that for any  $g:T\to S$ the following are equivalent
\begin{enumerate}
\item  
$\omega_{X\times_ST/T}^{[m]}$
is flat over $T$ and commutes with base change.
\item  $g$ factors through $S_m\to S$.
\end{enumerate}
\end{cor}

Proof. 
We claim that $S_m={\rm Hull}\bigl(\omega_{X/S}^{\otimes m}\bigr)$.
Given $g:T\to S$, let $j_T:X\times_ST\setminus Z\times_ST\to X\times_ST$
be the inclusion. Then 
$$
\omega_{X\times_ST/T}^{[m]}=\bigl(j_T\bigr)_*
g_X^*\omega_{X\setminus Z/S}^{\otimes m}.
$$
If $T\mapsto \omega_{X\times_ST/T}^{[m]}$
 commutes with restrictions to
the fibers of $X\times_ST\to T$, then 
$\omega_{X\times_ST/T}^{[m]}$ has $S_2$ fibers, hence
$\omega_{X\times_ST/T}^{[m]}$
 is the hull of $\omega_{X\times_ST/T}^{\otimes m}$.

Conversely, 
if $\omega_{X\times_ST/T}^{\otimes m}$ has a hull
then it is $\omega_{X\times_ST/T}^{[m]}$ and
it commutes with further base changes by (\ref{hull.defn}).\qed

\medskip

In order to prove (\ref{omega.all.comm.cor}), choose $N$ such that
$\omega_{X_s}^{[N]}$ is locally free for every $s\in S$.
For $1\leq i\leq N$,
let $S_i\to S$ be as in (\ref{omega.m.comm.cor}).
Take $T$ to be the fiber product 
$S_1\times_S \cdots\times_S S_N\to S$.
Then $\omega_{X\times_ST/T}^{[m]}$
is flat over $T$ and commutes with base change for every $m$.
Thus $S_H$ is the disjoint union of those connected components of
$T$ where the Hilbert function is $H$.\qed

\section{Coarse moduli spaces}\label{sec.5}

Having defined the correct moduli functor for
slc models, we can now get down to studying its properties
and the corresponding moduli spaces.

\begin{say}[Valuative criterion of separatedness] 

The functor $SlcMod$ satisfies the valuative criterion
of separatedness, which is essentially (\ref{ample-mm.prop}). This was built
into our construction.
\end{say}

\begin{say}[Valuative criterion of  properness] \label{5.proper}

The short answer is that  $SlcMod$ satisfies the valuative criterion
of properness, but some warnings are in order.

We proceed very much as for curves. We start with a family
of canonical models over an open  curve $X^0\to B^0$. 
By the semi-stable reduction theorem of \cite{kkms},
after a base change $C^0\to B^0$ and extending the family
over a proper curve $C\supset C^0$, there is a resolution
$g:Y\to C$ all of whose fibers are reduced simple normal crossing divisors.
Finally we replace $Y$ by its relative canonical model
$$
Y^c:=\proj_C \sum_{m\geq 0}g_* \bigl(\omega_{Y/C}^m\bigr).
$$
It is not hard to see that $Y^c\to C$ 
 is  in $SlcMod(C)$ extending $X\times_{B^0}C^0\to C^0$.

This establishes the valuative criterion of  properness
if canonical models are dense in the moduli of slc models.
We probably mostly care about the irreducible components
where canonical models are dense, so we could take this as
the final answer.

However, not all irreducible components are such, and it would be better to
understand all of them. 

So let us start with a family of slc models $X^0\to B^0$.
We can proceed as above, but instead of the relative canonical model
of $Y$ we need to take the relative semi log canonical model.
As we noted,  semi log canonical rings
 are not always finitely generated \cite{k-ncsurf}.

Here the solution is to normalize the family, construct the
models of the normalization over $C$ and then try to reconstruct
the desired extension of the original family.
This is actually quite subtle, see \cite[Chap.3]{k-m-book}.
\end{say}

\begin{say}[Existence of coarse moduli spaces]

Fix a  function $H$ and an integer $m$.
Let $SlcMod_{H,m}$ be the moduli functor of
 slc models with Hilbert funcion $H$ for which $\omega^{[m]}$ is
locally free, very ample and has no higher cohomologies.
 All of these thus embedd into
$\p^N$ for $N=H(m)-1$.  We use a variant of
(\ref{omega.all.comm.cor}) to show that there is a
locally closed subscheme $S_{H,m}$ of the Hilbert scheme $\hilb(\p^N)$
that parametrizes families of $m$-canonically embedded
slc models with Hilbert function $H$.  

The general quotient theorems of \cite{k-quot}, \cite{ke-mo}
apply and we obtain the coarse moduli space ${\rm SlcMod}_{H,m}$
of $SlcMod_{H,m}$ as the geometric quotient $S_{H,m}/\aut(\p^N)$.

Finally we let $m$ run through the sequence $2!, 3!, 4!, \dots$
to get an increasing sequence of coarse moduli spaces
whose union is the coarse moduli space ${\rm SlcMod}_{H}$.
 For now we know only that it is a
separated  algebraic space
which is locally of finite type.
\end{say}

\begin{say}[Properness]  We saw that ${\rm SlcMod}_{H}$
satisfies the valuative criterion of properness, hence it is proper iff
it is of finite type. 

The components where the canonical models are dense were
studied by \cite{karu}. He proves that one can control the
procedure outlined in (\ref{5.proper}) uniformly. Thus every such component is
of finite type, hence proper.

With some modifications, this implies that every irreducible
component of ${\rm SlcMod}_{H}$ is proper.

Thus the only remaining question is: can there by infinitely many
irreducibe components?

To illustrate some of the difficulties, let us consider
a much simpler question: can we bound the number of irreducible
components of a slc surface $S$ with Hilbert function $H$?

For curves the answer is easy. If $C=\cup_i C_i$
then
$$
2g(C)-2=\deg \omega_C=\tsum_i \deg \bigl(\omega_C|_{C_i}\bigr).
$$
Each $C_i$ on the right hand side contributes at least 1,
hence there are at most $2g-2$ irreducible
components. In the surface case, we have something very similar.
If $S=\cup_i S_i$ then we can compute the self intersection of the 
canonical class as
$$
\bigl(K_S\cdot K_S\bigr)=
\tsum_i \bigl(K_S|_{S_i}\cdot K_S|_{S_i}\bigr).
$$
The unexpected propblem is that $K_S$ is only a $\q$-Cartier divisor,
hence each summand on the right hand side is a
positive rational number, not an integer.

We are, however, saved if the contributions on the right
are bounded away from 0. This, and much more that is needed for
boundedness was proved in \cite{ale-bound} and improved in
\cite{al-mo}.
The lower bound  $\frac1{1764}$ was established in \cite{k-logsurf}.
(I do not know the optimal bound,
 but $\bigl(\p(3,4,5), (x^3y+y^2z+z^2x=0)\bigr)$
has $(K_S+D)^2=\frac1{60}$.)

In higher dimensions, recent work of Hacon-M\textsuperscript{c}Kernan-Xu
establishes a lower bound for $\bigl(K_X|_{X_i}\bigr)^n$;
the methods are likely to give boundedness as well.
\end{say}

\begin{say}[Projectivity] The method of \cite{k-proj} shows that
every proper subscheme of ${\rm SlcMod}$ is projective.
For $m$ sufficiently divisible, the 1-dimensional vector spaces
$\det H^0\bigl(X, \omega_X^{[m]}\bigr)$ naturally glue together
to an ample line bundle.

The proof uses the Nakai-Moishezon ampleness criterion,
thus it works only for proper schemes.

It seems very hard to give quasi-projectivity criteria.
For instance, \cite{k-ncsurf} gives an example of a normal crossing surface $S$
with a line bundle $L$ and normalization $\pi:\bar S\to S$ such that
$\pi^*L$ is ample yet $L$ is not ample, in fact 
no power of $L$ is generated by global sections.
\end{say}

\section{Moduli  of slc pairs}\label{sec.6}

In dimension 1, it is useful to consider not just
the moduli of curves but also the moduli of pointed curves.
Similarly, in higher dimensions, one should consider
the moduli of pairs $(X,\Delta)$ where $\Delta=\sum_i a_iD_i$ is a
linear combination of divisors with coefficients $0\leq a_i\leq 1$. 
These were first considered in \cite{ale-pairs}.

The first task is to define slc singularities of pairs. This is actually
quite natural, see \cite{k-etal}.

By contrast, 
finding the correct analog of (\ref{can.mod.funct.hilb.2.def})
 turns out to be a quite thorny
problem. Instead of going into  details, let me just present
a key example, due to Hassett, which shows
that in general we can not view a deformation of a pair
as first a deformation of $X$ and then a deformation of the $D_i$.
One must view $(X,\Delta)$ as an inseparable unit.

\begin{exmp}\label{nonflat.exmps}
Let $S\subset \p^5$ be the 
 cone  over the degree 4 rational normal curve.
Fix $r\geq 1$ and let $D_S$ be the sum of $2r$ lines.
 Then $(S,\frac1{r}D_S)$ is lc and
$\bigl(K_S+\frac1{r}D_S\bigr)^2=4$. 

There are two different deformations of the pair $(S,D_S)$.

(\ref{nonflat.exmps}.1) First, set $P:=\p^2$ and 
let $D_P$ be the sum of $r$ general lines. 
Then $(P,\frac1{r}D_P)$ is lc (even canonical if $r\geq 2$) and
$\bigl(K_P+\frac1{r}D_P\bigr)^2=4$. 
The usual smoothing of $S\subset \p^2$ to the Veronese 
surface gives a family $f:(X,D_X)\to \p^1$ with general fiber  $(P,D_P)$
and special fiber $(S,D_S)$. We can concretely realize this as
deforming $(P,D_P)\subset \p^5$ to the cone over a general
hyperplane section.
Note that for any general $D_S$ there is a choice of lines
$D_P$ such that the above limit is exactly $D_S$.

The total space  $(X,D_X)$ is the cone over
$(P,D_P)$ (blown up along  curve)
and $X$ is $\q$-factorial. 
The structure sheaf of an effective divisor on $X$ is CM.

In particular,
 $D_S$ is a flat limit of $D_P$. 
Since the $D_P$ is a plane curve of degree $r$,
we conclude that 
$$
\chi(\o_{D_S})=\chi(\o_{D_P})=-\frac{r(r-3)}{2}.
$$

(\ref{nonflat.exmps}.2) Second, 
 set $Q:=\p^1\times \p^1$ and
 let $A,B$ denote the classes of the 2 rulings.
 Let $D_Q$ be the sum of $r$ lines from the $A$-family.
 Then $(Q,\frac1{r}D_Q)$ is  canonical  and
$\bigl(K_Q+\frac1{r}D_Q\bigr)^2=4$. The usual smoothing of $S\subset \p^2$ to 
$\p^1\times \p^1$ embedded by $H:=A+2B$
 gives a family $g:(Y,D_Y)\to \p^1$ with general fiber  $(Q,D_Q)$
and special fiber $(S,D_S)$.
We can concretely realize this as
deforming $(Q,D_Q)\subset \p^5$ to the cone over a general
hyperplane section.

The total space $(Y,D_Y)$ is the cone over
$(Q,D_Q)$ (blown up along  curve) and $Y$ is not $\q$-factorial.
However, $K_Q+\frac1{r}D_Q\simq -H$, thus $K_Y+\frac1{r}D_Y$
 is $\q$-Cartier and
$(Y, \frac1{r}D_Y)$ is lc.

In this case, however, $D_Q$ is not a flat limit of $D_P$
for $r>1$. Thus follows, for instance, from
comparing their Euler characteristic:
$$
 \chi(\o_{D_S})=-\frac{r(r-3)}{2}\qtq{and}
 \chi(\o_{D_Q})=r.
$$

(\ref{nonflat.exmps}.3) Because of their role in the
canonical algebra, we are
also interested in the sheaves
$\o(mK+\rdown{\frac{m}{r}D})$. 

Let $H_P$ be the hyperplane class  of $P\subset \p^5$
(that is, 2 times a line  $L\subset P$) and
write $m=br+a$ where $0\leq a <r$.  
One computes that
$$
\begin{array}{ccl}
\chi\bigl(P,\o_P(mK_P+\rdown{\tfrac{m}{r}D_P}+nH_P)\bigr)&=&
\binom{2n-2m+2}{2} -a(2n-2m+1)+\binom{a}{2},\\[1ex]
\chi\bigl(S,\o_S(mK_S+\rdown{\tfrac{m}{r}D_S}+nH_S)\bigr)&=&
\binom{2n-2m+2}{2} -a(2n-2m+1)+\binom{a}{2},\\[1ex]
\chi\bigl(Q,\o_Q(mK_Q+\rdown{\tfrac{m}{r}D_Q}+nH_Q)\bigr)&=&
\binom{2n-2m+2}{2} -a(2n-2m+1).
\end{array}
$$
From this we conclude that  the restriction of $\o_Y(mK_Y+\rdown{mD_Y})$
to the central fiber $S$ agrees with 
$\o_S(mK_S+\rdown{mD_S})$  only if $a\in \{0,1\}$,
that is when $m\equiv 0,1\mod r$. 
The if part was clear from the beginning. Indeed, if $a=0$ then 
$\o_Y(mK_Y+\rdown{mD_Y})$ is locally free and if
 $a=1$ then 
$\o_Y(mK_Y+\rdown{mD_Y})$ is $\o_Y(K_Y)$
tensored with a locally free sheaf. Both of these
commute with restrictions.

In the other cases we only get an injection
$$
\o_Y(mK_Y+\rdown{mD_Y})|_S\into \o_S(mK_S+\rdown{mD_S})
$$
whose quotient is a torsion sheaf of length $\binom{a}{2}$
supported at the vertex.
\end{exmp}

There are several ways to overcome these problems;
all them will be discussed in \cite{k-m-book}.

\begin{enumerate}
\item  Embedded points do not appear if all the coeffcients $a_i$
are $>\frac12$. 
\item By wiggling the $a_i$ suitably, one again avoids
embedded points.
\item  Fix $m$ such that $\o_X(mK_X+m\Delta)$ is locally free.
One can identify a pair $(X,\Delta)$
with the corresponding map $\omega_X^{\otimes m}\to \o_X(mK_X+m\Delta)$.
It turns out to be easier to  deal with the moduli of triples
$\bigl(X, \omega_X^{\otimes m}\to L\bigr)$ for some line bundle $L$.
\item The branch varieties of \cite{ale-knu} give another approach.
\end{enumerate}

 \begin{ack} Parts of this paper were written in connection with a 
lecture series on moduli at IHP, Paris.
I thank my audience and especially C.~Voisin for the invitation,
 useful comments and corrections.
Partial financial support  was provided by  the NSF under grant number 
DMS-0758275.
\end{ack}

\bibliography{refs}

\vskip1cm

\noindent Princeton University, Princeton NJ 08544-1000

\begin{verbatim}kollar@math.princeton.edu\end{verbatim}

\end{document}